\documentclass[psamsfonts, reqno, 12pt, letterpaper]{amsart}
\usepackage{amsfonts,amsmath, amsthm, amssymb, latexsym, epsfig}
\usepackage[all]{xy}
\usepackage{xspace}
\usepackage{comment}
\usepackage{setspace}
\usepackage{enumerate}
\usepackage[mathscr]{eucal}

	\advance \topmargin by -\headheight
	\advance \topmargin by -\headsep

	\evensidemargin \oddsidemargin
	\newcommand{\Ntau}{\mathsf{n}_\tau}
	\newcommand{\Ntauo}{\mathsf{n}_{\tau_1}}
	\newcommand{\Nupsilon}{\mathsf{n}_{\tau_2}}
	\newcommand{\ppvec}{\mathsf{p}}

	\newcommand{\ftn}[3]{ #1 : #2 \rightarrow #3 }

	\newcommand{\norm}[1]{\left\|#1\right\|}

	\newcommand{\kthy}{\ensuremath{\mathnormal{K}}-theory\xspace}

	\newcommand{\unit}[1]{1_{ #1 }}
	
	\newcommand{\ksix}{\ensuremath{K_{\mathrm{six}}}\xspace}
		\newcommand{\Ksix}{\ensuremath{\mathbf{K}_{\mathrm{six}}}\xspace}

	\newcommand{\mc}[1]{\mathcal{#1}}
	
	\newcommand{\msf}[1]{\mathsf{#1}}

	\newcommand{\mfk}[1]{\mathfrak{#1}}
	\newcommand{\mbf}[1]{\mathbf{#1}}
	
	\newcommand{\e}{\ensuremath{\mathfrak{e}}\xspace}
	\newcommand{\Z}{\ensuremath{\mathbb{Z}}\xspace}
	\newcommand{\C}{\ensuremath{\mathbb{C}}\xspace}
	\newcommand{\Q}{\ensuremath{\mathbb{Q}}\xspace}
	
	\newcommand{\N}{\ensuremath{\mathbb{N}}\xspace}
	
	\newcommand{\K}{\ensuremath{\mathcal{K}}\xspace}

	\newcommand{\kk}{\ensuremath{\mathit{KK}}\xspace}

	\newcommand{\coker}{\ensuremath{\operatorname{coker}}}

	\newcommand{\Hom}{\ensuremath{\operatorname{Hom}}}

	\newcommand{\id}{\ensuremath{\operatorname{id}}}
	
	\newcommand{\bootstrap}{\ensuremath{\mathcal{N}}\xspace}
	\newcommand{\Ell}{K_*^+}
        	\newcommand{\Extab}{\ensuremath{\operatorname{Ext}}_{\Z}^{1}}
	
	\newcommand{\sixs}{\ensuremath{\mathrm{Hext}}}
	
	\newcommand{\Exts}{\ensuremath{\mathcal{E}\mathrm{xt}}}
	\newcommand{\Ext}{\ensuremath{\operatorname{Ext}}}

	\newcommand{\multialg}[1]{\mathcal{M}(#1)\xspace}
	\newcommand{\corona}[1]{\mathcal{Q}(#1)\xspace}
	\newcommand{\starhom}{{$*$}-ho\-mo\-mor\-phism\xspace}
	
	\newcommand{\stariso}{{$*$}-iso\-mor\-phism\xspace}
	\newcommand{\starisos}{{$*$}-isomorphisms\xspace}

		\newcommand{\image}{\ensuremath{\operatorname{image}}\xspace}

	\newcommand{\cstar}{{$C \sp \ast$}\xspace}
	\newcommand{\TSS}[1][]{{{\underline{\mathsf{X}}_{#1}}}}
	\newcommand{\mattau}[1][]{\mathcal{O}_{\TSS[#1]}}

	\theoremstyle{plain}
	\newtheorem{thm}{Theorem}[section]
	\newtheorem{lemma}[thm]{Lemma}
	\newtheorem{theorem}[thm]{Theorem}
	\newtheorem{proposition}[thm]{Proposition}
	\newtheorem{corollary}[thm]{Corollary}

	\newtheorem{assumption}[thm]{Assumption}
	
	\theoremstyle{definition}
	\newtheorem{definition}[thm]{Definition}
	
	\newtheorem{remark}[thm]{Remark}
	
	\newtheorem{example}[thm]{Example}

	\numberwithin{equation}{section}
	\numberwithin{figure}{section}
\newenvironment{skproof}{{\noindent\emph{Sketch of proof.}}}{\hfill$\square$}

\begin{document}
	\title[Extensions of Classifiable \cstar-algebras]{Classification of Extensions of Classifiable \cstar-algebras}
	\author{S{\o}ren Eilers}
        \address{Department of Mathematical Sciences \\
        University of Copenhagen\\
        Universitetsparken~5 \\
        DK-2100 Copenhagen, Denmark}
        \email{eilers@math.ku.dk }
        \author{Gunnar Restorff}
        \address{Department of Mathematical Sciences \\
        University of Copenhagen\\
        Universitetsparken~5 \\
        DK-2100 Copenhagen, Denmark}
\curraddr{Faculty of Science and Technology\\University of Faroe 
Islands\\N\'oat\'un 3\\FO-100 T\'orshavn\\Faroe Islands}
\email{gunnarr@setur.fo}

	\author{Efren Ruiz}
        \address{Department of Mathematics\\University of Hawaii,
Hilo\\200 W. Kawili St.\\
Hilo, Hawaii\\
96720-4091 USA}
        \email{ruize@hawaii.edu}
        \date{Third revised version, \today}
	


	\begin{abstract}
	For a certain class of extensions $\e : 0 \to B \to E \to A
\to 0$ of \cstar-algebras in which $B$ and $A$ belong to 
classifiable classes of \cstar-algebras, we show that the functor which
sends $\e$ to its associated six term exact sequence in \kthy and the
positive cones of $K_{0}(B)$ and $K_{0} (A)$ is a classification
functor.  We give two independent applications addressing the
classification  of a class of
\cstar-algebras arising from substitutional shift spaces  on one hand and
 of graph
algebras
on the
other. The former application leads to the answer of a question of
Carlsen and the
first named author concerning the completeness of stabilized Matsumoto
algebras as an invariant of flow equivalence. The latter leads to the first classification result for non-simple 
graph $C*$-algebras.
	\end{abstract}
        \maketitle
        
 \section*{Introduction}

The magnificent recent progress of the classification theory for
simple $C^*$-algebras has few direct consequences for general
$C^*$-algebras, even for
those with finite ideal lattices. Furthermore, it is not even clear
what kind of $K$-theoretical invariant to use in such a context.

When there is just one non-trivial ideal, however, there is a
canonical choice of invariant.
Associated to every extension $0 \to B \to E \to A \to 0$ of nonzero
\cstar-algebras is the standard six term exact sequence of $K$-groups 
\begin{equation*}
\xymatrix{
K_{0} ( B ) \ar[r] & K_{0} ( E ) \ar[r] & K_{0} ( A ) \ar[d] \\
K_{1} ( A ) \ar[u] & K_{1} ( E ) \ar[l] & K_{1} ( B ) \ar[l],
}
\end{equation*}
providing a necessary condition for two extensions to be isomorphic. For examples of classification results involving the six term exact sequence of $K$-groups see \cite{gentoe}, \cite{inftoeplitz}, \cite{LinK1zero}, and \cite{inftoeplitzE}.  In each case, the extensions considered were extensions that can be expressed as inductive limits of simpler extensions.  The classification results were achieved by using the standard intertwining argument.    

In \cite{extpurelyinf}, R{\o}rdam used a completely different
technique to classify a certain class of extensions.  He considered
essential extensions of separable nuclear purely infinite simple
\cstar-algebras in \bootstrap, where \bootstrap is the bootstrap
category of Rosenberg and Schochet \cite{uct}.  Employing the fact that
every invertible element of $\kk ( A , B )$ (where $A$ and $B$ are
separable nuclear stable purely infinite simple \cstar-algebras) lifts
to a \stariso from $A$ to $B$ and that every essential extension of
$A$ by $B$ is absorbing, R{\o}rdam showed that the six term sequence
is, indeed, a complete invariant in this case.

The purpose of this paper is to extend the above result to other
classes of \cstar-algebras that are classified via $K$-theoretical
invariants. As we shall see, both the celebrated classification results of
Kirchberg-Phillips (\cite{ekncp:eeccao}) and of Lin (\cite{hl:cscttrz}) can be transfered to this setting under
an assumption of \emph{fullness} of the extension which is automatic
in the case solved by R{\o}rdam.
 
The motivation of our work was an application to a class of
\cstar-algebras introduced in the work of Matsumoto. In a case studied
in \cite{CarlSym}, \cite{CarlEilersKgrpsMat} one gets that the Matsumoto algebra $\mattau[]$ fits in a short
exact sequence of the form
\begin{equation*}
\xymatrix{ 0 \ar[r] & \K^n \ar[r] & \mattau[] \ar[r] & C(\TSS) \rtimes_{\sigma} \Z \ar[r] & 0. }
\end{equation*}
Since $C(\TSS) \rtimes_{ \sigma } \Z$ is a unital simple $AT$-algebra
with real rank zero our results apply to classify $\mattau[]$ by it's
$K_0$-group with a scale consisting of $n$ preferred elements.    
 
The paper is organized as follows.  In Section \ref{extensions}, we
give basic properties and develop some notation concerning extensions
of \cstar-algebras.  Section \ref{sixk} gives notation (mainly from
\cite{extpurelyinf}) concerning the
six term exact sequence of $K$-groups and extends work of R{\o}rdam.
  Section \ref{classresult} contains
our main results (Theorem \ref{T : classification} and Theorem \ref{T
: matsualgs}).  In the last section we use these results to classify
the \cstar-algebras described in the previous paragraph. We also
present an alternative application to \emph{graph $C^*$-algebras} which fully employs the capacity of our classification result to handle $C^*$-algebras which have some subquotients which are stably finite, and some which are purely infinite.

An earlier
  version of this paper was included in the second named author's
  PhD thesis.
 
 \section{Extensions}\label{extensions}

\subsection{Notation}

    For a stable $C\sp*$-algebra $B$ and a $C\sp*$-algebra $A$, we
    will denote the class of essential extensions 
    $$0\rightarrow B\overset{\varphi}{\rightarrow} E\overset{\psi}{\rightarrow} A\rightarrow 0$$
    by $\mathcal{E}\mathrm{xt}(A, B)$.

    Since the goal of this paper is to classify extensions of
    separable nuclear $C\sp*$-algebras, throughout the rest of the
    paper we will only consider $C\sp*$-algebras that are separable
    and nuclear.  
    
    \begin{assumption}
   {In the rest of the paper all 
    $C\sp*$-algebras considered are assumed to be separable and
    nuclear unless stated otherwise. Note in particular that 
    multiplier and corona algebras will be non-separable.}  
\end{assumption}

    Under the above assumption, if $B$ is a stable $C\sp*$-algebra, then
    we may identify $\mathrm{Ext}(A, B)$ with $KK^1(A,B)$ 
    (for the definition of $\mathrm{Ext}(A,B)$ and $KK^i(A,B)$
    see Chapters 7 and 8 in \cite{blackadarB}). 
    So for $x$ in $\mathrm{Ext}(A,B)$ and $y$ in $KK^i(B,C)$, the
    Kasparov product $x\times y$ is an element of $KK^{i+1}(A,C)$. 
    For every element $\mathfrak{e}$ of $\mathcal{E}\mathrm{xt}(A,B)$,
    we use $x_{A,B}(\mathfrak{e})$ to denote the element of
    $\mathrm{Ext}(A,B)$ that is represented by $\mathfrak{e}$. \\

 \begin{definition}\label{D : morphism1}
    A homomorphism from an extension 
    $0\rightarrow B_1\rightarrow E_1\rightarrow A_1\rightarrow 0$ 
    to an extension 
    $0\rightarrow B_2\rightarrow E_2\rightarrow A_2\rightarrow 0$ 
    is a triple $(\beta,\eta,\alpha)$ such that the diagram
    $$\xymatrix{
      0\ar[r]&B_1\ar[r]\ar[d]_{\beta}&E_1\ar[r]\ar[d]_{\eta}&A_1\ar[r]\ar[d]_{\alpha}&0 \\
      0\ar[r]&B_2\ar[r]&E_2\ar[r]&A_2\ar[r]&0 
    }$$
    commutes. 
  \end{definition}
  
    This turns the class of extensions of $C\sp*$-algebras into a
    category in the canonical way.

    We say that an extension 
    $\mathfrak{e}_1\colon 
    0\rightarrow B\rightarrow E_1\rightarrow A\rightarrow 0$ 
    is \emph{congruent} to an extension
    $\mathfrak{e}_2\colon 
    0\rightarrow B\rightarrow E_2\rightarrow A\rightarrow 0$, 
    if there exists an isomorphism of the form
    $(\mathrm{id}_B,\eta,\mathrm{id}_A)$ from $\mathfrak{e}_1$ to $\mathfrak{e}_2$.

We will use the following notation from \cite{extpurelyinf}. For each injective \starhom $ \alpha$ from $A_{1} $ to $A_{2} $ and for each $\e$ in $\Exts( A_{2}, B)$, there exists a unique extension $\alpha \cdot \e$ in $\Exts ( A_{1}, B )$ such that the diagram
\begin{equation*}
\xymatrix{ 
\alpha \cdot \e : & 0 \ar[r] & B \ar[r] \ar@{=}[d] & \widetilde{E} \ar[r] \ar@{^{(}->}[d] & A_{1} \ar[r] \ar[d]_{\alpha} & 0\\
				\e : & 0 \ar[r] & B \ar[r]  & E \ar[r]  & A_{2} \ar[r]  & 0 }
\end{equation*}
is commutative.  For each \stariso $ \beta $ from $B_{1}$ to $B_{2}$ and for each $\e$ in $\Exts ( A, B_{1} )$, there exists a unique extension $\e \cdot \beta$ in $\Exts ( A, B_{2} )$ such that the diagram
\begin{equation*}
\xymatrix{ 
 \e			 : & 0 \ar[r] & B_{1} \ar[r] \ar[d]_{\beta} & E \ar[r] \ar@{=}[d] & A \ar[r] \ar@{=}[d] & 0\\
\e \cdot \beta	 : & 0 \ar[r] & B_{2} \ar[r]  & E \ar[r]  & A \ar[r]  & 0
				 }
\end{equation*}     
is commutative.

Proposition 1.1 and Proposition 1.2 in \cite{extpurelyinf} explain
the interrelations between the concepts introduced above and will be
crucial in our approach.

\subsection{Full extensions}\label{fullextchap}

Let $a$ be an element of a \cstar-algebra $A$.  We say that \emph{$a$
is norm-full in $A$} if $a$ is not contained in any norm-closed proper
ideal of $A$.  The word ``full'' is also widely used, but since we
will often work in multiplier algebras, we emphasize that it is the
norm topology we are using, rather than the strict topology.  The next
lemma is a consequence of a result of L.G.\ Brown (see Corollary 2.6
in \cite{heralgs}); we leave the proof to the reader.

\begin{lemma}\label{L : brown}
Let $A$ be a separable \cstar-algebra.  If $p$ is a norm-full projection in $A \otimes \mbf{M}_{n} \subset A \otimes \K$, then there exists a \stariso $\varphi$ from $A \otimes \K$ onto $p ( A \otimes \K ) p \otimes \K$ such that $[ \varphi ( p ) ] = [ p \otimes e_{11} ]$.   
\end{lemma}

\begin{definition}\label{fullext}
An extension $\e$ is said to be \emph{full} if the associated Busby
invariant $\tau_\e$ has the property that $\tau_\e(a)$ is a norm-full
element of $\corona{B}$ for any $a\in A\backslash\{0\}$.
\end{definition}

\begin{lemma}\label{littlel} If $B$ is stable and purely infinite, then any extension
in $\Exts(A,B)$ is full. If $B$ is stable and $A$ is unital and simple, then any
unital extension in $\Exts(A,B)$ is full.
\end{lemma}
\begin{proof}
In the first case, the corona algebra is simple
(\cite{simplecorona}). In the second, the image of the Busby map is a simple
unital subalgebra of the corona algebra  and hence
cannot intersect an ideal nontrivially.
\end{proof}

It seems reasonable to expect that  the stabilized extension
of a full extension is again full. We prove this under the added
assumption that $B$ is already stable:

\begin{proposition}\label{P : stablefull}
Let $\e :\ 0 \to B \overset{ \iota }{ \to } E \overset{ \pi }{ \to } A
\to 0$ be an essential extension where $B$ is a stable
\cstar-algebra. If $\e$ is full, then so is
\begin{equation*}
\xymatrix{
\e^{s} : & 0 \ar[r] & B \otimes \K \ar[r]^{ \iota \otimes \id_{ \K } } & E \otimes \K \ar[r]^{ \pi \otimes \id_{ \K } } &  A \otimes \K \ar[r] & 0}
\end{equation*}  
\end{proposition}

\begin{proof}
For any \cstar-algebra $C$, denote the embedding of $C$ into $C
\otimes \K$ which sends $c$ into $c \otimes e_{11}$ by $\iota_{C}$ and denote the canonical embedding of $C$ as an essential ideal of the multiplier algebra $\multialg{C}$ of $C$ by $\theta_{C}$.  We will first show that $\iota_{B}$ satisfies the following properties:
\begin{enumerate}
\item $\iota_{B}$ has an extension $\widetilde{ \iota }_{B}$ from $\multialg{B}$ to $\multialg{ B \otimes \K }$ (i.e.\ $\theta_{ B \otimes \K } \circ \iota_{ B } = \widetilde{ \iota }_{B}  \circ \theta_{ B } $), which maps $1_{ \multialg{ B } }$ to a norm-full projection in $\multialg{ B \otimes \K }$ and

\item the map $\bar{ \iota }_{ B }$ from $\corona{B}$ to $\corona{ B \otimes \K }$ induced by $\tilde{ \iota }_{B}$ intertwines the Busby invariants of $\e$ and $\e^{s}$, and the $*$-homomorphism $\iota_{ A }$ (i.e.\ $\tau_{ \e^{s} } \circ \iota_{A} = \overline{ \iota }_{ B } \circ \tau_{ \e }$).   
\end{enumerate}

First note that there exist unique injective $*$-homomorphisms $\sigma$ from $E$ to $\multialg{B}$ and $\sigma^{s}$ from $E \otimes \K$ to $\multialg{ B \otimes \K}$ such that $\theta_{B} = \sigma \circ \iota$ and $\theta_{ B \otimes \K } = \sigma^{s} \circ ( \iota \otimes \id_{ \K } )$.  It is well-known that we have a unique $*$-homomorphism $\rho$ from $\multialg{B} \otimes \multialg{ \K }$ to $\multialg{B \otimes \K }$ such that $\theta_{ B \otimes \K } = \rho \circ ( \theta_{B} \otimes \theta_{ \K } )$ and that this map is injective and unital (see Lemma 11.12 in \cite{PedPullPush}).

In the following diagram, all the maps are injective $*$-homomorphisms
\begin{equation*}
\xymatrix{
		& & B \otimes \K \ar@/_/[lldd]_-{ \theta_{B} \otimes \id_{ \K } } \ar@/^/[rrdd]^-{ \theta_{ B \otimes \K } } \ar[d]_-{ \iota \otimes \id_{ \K } } &  &   \\
	& & E \otimes \K \ar[lld]^-{ \sigma \otimes \id_{ \K } } \ar[rrd]_-{ \sigma^{s} } & &   \\
\multialg{B} \otimes \K \ar[rr]_-{ \id_{ \multialg{B} } \otimes \theta_{\K} } & & \multialg{B} \otimes \multialg{K} \ar[rr]_-{ \rho } & &\multialg{ B \otimes \K }
}
\end{equation*}
The bottom triangle commutes  by the uniqueness of $\sigma^{s}$, so
this is a commutative diagram.

Now let $\widetilde{ \iota }_{ B } = \rho \circ ( \id_{ \multialg{B} } \otimes \theta_{ \K } ) \circ \iota_{ \multialg{B} }$.  Clearly, $\theta_{ B \otimes \K } \circ \iota_{B} = \widetilde{ \iota }_{ B } \circ \theta_{B}$ and $p = \widetilde{ \iota }_{ B } ( 1_{ \multialg{B} } )$ is a projection in $\multialg{ B \otimes \K }$.  Note that $\iota_{B} ( B ) = B \otimes e_{11} \cong p \theta_{ B \otimes \K } ( B \otimes \K ) p$.  Therefore, $p \theta_{ B \otimes \K } ( B \otimes \K ) p$ is a stable, hereditary, sub-\cstar-algebra of $\theta_{ B \otimes  \K } ( B \otimes \K )$ which is not contained in any proper ideal of $\theta_{ B \otimes \K } ( B \otimes \K )$.  By Theorem 4.23 in \cite{BrownSemi}, $p$ is Murray-von Neumann equivalent to $1_{ \multialg{ B \otimes \K } }$.  Hence, $p  = \widetilde{ \iota }_{ B } ( 1_{ \multialg{ B } } )$ is norm-full in $\multialg{ B \otimes \K }$.  

Now we see that $\widetilde{ \iota }_{ B } \circ \sigma = \sigma^{s} \circ \iota_{ E }$ since the following diagram is commutative:
\begin{equation*}
\xymatrix{
E \ar[d]_{ \sigma } 	\ar[rrrr]^-{ \iota_{E} }				& 														& & & E \otimes \K \ar[d]^-{ \sigma^{s} } \ar[dlll]_-{ \sigma \otimes \id_{ \K } } \\
\multialg{ B } \ar[r]_-{ \iota_{ \multialg{B} } } & \multialg{B} \otimes \K \ar[rr]_-{ \id_{ \multialg{B} } \otimes \theta_{\K } } & &\multialg{B} \otimes \multialg{ \K } \ar[r]_-{ \rho } & \multialg{ B \otimes \K } 
}
\end{equation*}      
  
Let $\overline{ \iota }_{ B }$ denote the $*$-homomorphism from $\corona{B}$ to $\corona{ B \otimes \K }$ which is induced by $\widetilde{ \iota }_{ B }$.  Arguing as in the proof of Theorem 2.2 in \cite{ELPmorph}, we have that the diagram
\begin{equation}\label{diagram}
\vcenter{\xymatrix{
A \ar[r]^{ \tau_{ \e } } \ar[d]_{ \iota_{A} } 	& \corona{ B } \ar[d]^{ \overline{ \iota }_{ B } } \\
A \otimes \K \ar[r]_{ \tau_{ \e^{s} } } 		& \corona{ B \otimes \K }
}}
\end{equation}    
is commutative since $( \iota_{B} , \iota_{E} , \iota_{A} )$ is a morphism from $\e$ to $\e^{s}$.  This finishes the proof of the two claims (1) and (2) above.

We are now ready to prove the proposition.  Let $x$ be a nonzero positive element of $A \otimes \K$.  Then there exist $t$ and $s$ in $A \otimes \K$ such that $s x^{ \frac{1}{2} } t = \iota_{A} ( y )$ for some nonzero positive element $y$ of $A$.  Let $\epsilon$ be a strictly positive number.  From (1) of our claim, there exist $x_{1}, \dots, x_{n}$, $y_{1}, \dots , y_{n}$ in $\corona{ B \otimes \K }$ such that 
\begin{equation*}
\norm{  1_{ \corona{ B \otimes \K } } - \sum_{ i = 1 }^{n} x_{i} \overline{ \iota }_{B} ( 1_{ \corona{B} } ) y_{i} } < \frac{ \epsilon }{ 2 }.
\end{equation*}
From our assumption on $\tau_{ \e }$, there exist $t_{1},\dots, t_{m}$, $s_{1}, \dots, s_{m}$ in $\corona{ B }$ such that 
\begin{equation*}
\norm{ 1_{ \corona{B } }- \sum_{ j = 1 }^{m} s_{j} \tau_{ \e } ( y ) t_{j}  } < \frac{ \epsilon }{ 2 ( \sum_{ i = 1 }^{n} \norm{ x_{i} }\norm{ y_{i} } + 1 ) }.
\end{equation*}
An easy computation shows that
\begin{equation*}
\norm{ 1_{ \corona{ B \otimes \K } } - \sum_{ i = 1 }^{n} x_{i} \left( \sum_{j = 1 }^{m} \overline{\iota}_{B} ( s_{j} \tau_{ \e } ( y ) t_{j} ) \right) y_{i} } < \epsilon.
\end{equation*}

By the commutativity of Diagram (\ref{diagram}), we have that  
\begin{equation*}
\norm{ 1_{ \corona{ B \otimes \K } } - \sum_{ i = 1 }^{n} x_{i} \left( \sum_{j = 1 }^{m} \overline{\iota}_{B} ( s_{j} ) \tau_{ \e }^{s}( s x^{\frac{1}{2} } t )  \overline{\iota}_{B} ( t_{j} ) \right) y_{i} } < \epsilon.
\end{equation*}  
Therefore, the ideal of $\corona{ B \otimes \K }$ generated by $\tau_{ \e^{s} } ( x^\frac{1}{2} )$ is equal to $\corona{ B \otimes \K }$.  Since $x^{\frac{1}{2}}$ is contained in the ideal of $A \otimes \K $ generated by $x$, we have that $x$ is norm-full in $\corona{ B \otimes \K }$.  

For an arbitrary nonzero element $x$ of $A \otimes \K$, consider the positive nonzero element $x^{*} x$ of $A \otimes \K$ and apply the result on positive elements to conclude that $\tau_{ \e^{s} } ( x^{*} x )$ is norm-full in $\corona{ B \otimes \K }$.  Therefore, $\tau_{ \e^{s} } ( x )$ is norm-full in $\corona{ B \otimes \K }$ since $x^{*} x$ is contained in the ideal of $A \otimes \K$ generated by $x$.
\end{proof}

\section{Six term exact sequence in \kthy}\label{sixk}

    We need to extend some results by R\o rdam (\cite{extpurelyinf}) to a more general
    setting. 
    First we recall some of the definitions there.

    \subsection{} 
    For every $\mathfrak{e}$ in $\mathcal{E}\mathrm{xt}(A,B)$ we
    denote the cyclic six term exact sequence associated to
    $\mathfrak{e}$ by $K_{\mathrm{six}}(\mathfrak{e})$. 
    Let $\mathcal{H}\mathrm{ext}(A,B)$ denote the class of all cyclic
    six term exact sequences arising from elements of
    $\mathcal{E}\mathrm{xt}(A,B)$. 
    A homomorphism between such cyclic six term exact sequences is a
    6-periodic chain homomorphism. 
    We will frequently denote a homomorphism from 
    $K_{ \mathrm{six} } ( \mathfrak{e}_{1} )$, with $\mathfrak{e}_{1} : 0\rightarrow B_1\rightarrow E_1\rightarrow A_1\rightarrow 0$, to
    $K_{ \mathrm{six} } ( \mathfrak{e}_{2} )$, with
    $\mathfrak{e}_{2} : 0\rightarrow B_2\rightarrow E_2\rightarrow A_2\rightarrow 0$, by a
    triple $(\beta_*,\eta_*,\alpha_*)$, 
    where $\beta_*$ from $K_*(B_1)$ to $K_*(B_2)$, 
    $\eta_*$ from $K_*(E_1)$ to $K_*(E_2)$, and 
    $\alpha_*$ from $K_*(A_1)$ to $K_*(A_2)$ are homomorphisms
    (making the obvious diagrams commutative).

    We say that two elements $h_1$ and $h_2$ of
    $\mathcal{H}\mathrm{ext}(A,B)$ are \emph{congruent} if there is an
    isomorphism from $h_1$ to $h_2$ of the form
    $(\mathrm{id}_{K_*(B)},\eta_*,\mathrm{id}_{K_*(A)})$. 
    Let $\mathrm{Hext}(A,B)$ be the set of all congruence classes of
    $\mathcal{H}\mathrm{ext}(A,B)$. 
    For every element $h$ of $\mathcal{H}\mathrm{ext}(A,B)$ we let
    $\mathbf{x}_{A,B}(h)$ denote the congruence class in
    $\mathrm{Hext}(A,B)$ containing $h$. 
    According to \cite{extpurelyinf}, Proposition 2.1, there is a unique map
    $\mathbf{K}_{\mathrm{six}}$ from $\mathrm{Ext}(A,B)$ to
    $\mathrm{Hext}(A,B)$ such that
    $\mathbf{x}_{A,B}(K_{\mathrm{six}}(\mathfrak{e}))
    =\mathbf{K}_\mathrm{six}(x_{A,B}(\mathfrak{e}))$ for every element
    $\mathfrak{e}$ of $\mathcal{E}\mathrm{xt}(A,B)$.

    We let $K_*$ denote the map in the Universal Coefficient Theorem
    from $KK^i(A,B)$ to $\textrm{Hom}(K_*(A),K_{*+i}(B))$, for
    $i=0,1$. 
    Moreover, we set 
    \begin{align*}
      \mathrm{Hext}(A,B,\delta_*)
      &=\{h\in\mathrm{Hext}(A,B)\,|\,K_*(h)=\delta_*\},\\
      \mathrm{Ext}_{\delta_*}(A,B)
      &=\left\{x\in\mathrm{Ext}(A,B)\,\left|\,\begin{array}{c}\ker\delta_j\subseteq\ker(K_j(x))\\\mathrm{image}(K_j(x))\subseteq\mathrm{image}(\delta_j)
      \end{array}
      ,j=0,1\right\}\right. .
    \end{align*}
    \quad
    In the following, we will need the maps 
    $$(\sigma_{\delta_*}=)\sigma_{A,B,\delta_*}\colon
    \mathrm{Hext}(A,B;\delta_*)\rightarrow
    \mathrm{Ext}_\mathbb{Z}^1(\ker(\delta_*),\mathrm{coker}(\delta_{*+1}))$$
    $$(s_{\delta_*}=)s_{A,B,\delta_*}\colon
    \mathrm{Ext}_{\delta_*}(A,B)\rightarrow
    \mathrm{Ext}_\mathbb{Z}^1(\ker(\delta_*),\mathrm{coker}(\delta_{*+1}))$$
    introduced on pages 101--103 in \cite{extpurelyinf}.

\begin{lemma}\label{L : isoextdelta}
Let $A$, $B$, and $C$ be separable nuclear \cstar-algebras with $B$ stable and let $\delta_{*}$ be an element of $\Hom ( K_{*} ( C ) , K_{ * + 1 } ( B ) )$.  Suppose $C$ is in \bootstrap and suppose $x$ in $\kk ( A, C )$ is a $\kk$-equivalence.

Set $\lambda_{*} = \delta_{*} \circ K_{*} ( x )$ in $\Hom ( K_{*} ( A ) , K_{ * + 1 } ( B ) )$.  Then 
\begin{enumerate}
\item $x \times ( \cdot )$ is an isomorphism from $\Ext_{ \delta_{*} } ( C, B )$ onto $\Ext_{ \lambda_{*} } ( A , B )$.
\item $x$ induces an isomorphism $ [ K_{*} ( x ) ]$ from $\Extab ( \ker ( \delta_{*} ) , \coker ( \delta_{ * + 1 } ) )$ onto $\Extab ( \ker ( \lambda_{*} ) , \coker ( \lambda_{* + 1} ) )$.
\item Moreover, if $A$ and $B$ are in \bootstrap and if $x = \kk ( \alpha )$ for some injective \starhom $\alpha$ from $A$ to $C$, then the diagram
\begin{equation*}
\xymatrix{
\Ext_{ \delta_{*} } ( C, B )  \ar[rr]^{ x \times ( \cdot ) }_{\cong} \ar[d]_{ s_{ \delta_{*} } }	 & &\Ext_{ \lambda_{*} } ( A, B ) \ar[d]^{ s_{ \lambda_{*} } } \\
\Extab ( \ker ( \delta_{*} ) , \coker ( \delta_{*+1} ) ) \ar[rr]_{ [ K_{*} ( x ) ] }^{\cong}       & & \Extab ( \ker ( \lambda_{*} ) , \coker ( \lambda_{* + 1} ) )
}
\end{equation*}
is commutative.
\end{enumerate} 
\end{lemma}

\begin{proof}
Since $x$ is a $\kk$-equivalence, $x \times ( \cdot )$ is an isomorphism from $\Ext ( C, B )$ onto $\Ext ( A, B )$.  Therefore, to prove (1) it is enough to show that $x \times ( \cdot )$ maps $\Ext_{ \delta_{*} } ( C, B )$ to $ \Ext_{ \lambda_{*} } ( A, B )$ and $x^{-1} \times ( \cdot )$ maps $\Ext_{ \lambda_{*} } ( A, B )$ to $\Ext_{ \delta_{*} } ( C, B )$.  

Note that $K_{*} ( x )$ is an isomorphism and $K_{j} ( x \times z ) = K_{j} ( z ) \circ K_{j} ( x )$ for $j = 0 ,1$ and $z$ in $\Ext ( C, B )$.  Hence, $\image ( K_{j} ( z ) \circ K_{j} ( x ) ) = \image ( K_{j} ( z )  )$ and $\image ( \delta_{j} ) = \image ( \delta_{j} \circ K_{j} ( x ) )$.  By definition, if $z$ is in $\Ext_{ \delta_{*} } ( C, B )$, then $\image ( K_{j} ( z )  ) \subset  \image ( \delta_{j} )$ for $j = 0 , 1$.  Therefore, for $j = 0 ,1$, 
\begin{align*}
\image ( K_{j} ( z ) \circ K_{j} ( x ) )  \subset \image ( \delta_{j} \circ K_{j} ( x ) ) = \image ( \lambda_{j} ).
\end{align*}
A straightforward computation shows that $\ker ( \lambda_{j} ) \subset \ker ( K_{j} ( z ) \circ K_{j} ( x ) )$.  Hence, $x \times z$ is an element of $\Ext_{ \lambda_{*} } ( A, B )$ for all $z$ in $\Ext_{ \delta_{*} } ( C, B )$.  A similar computation shows that $x^{-1} \times ( \cdot )$ maps $\Ext_{ \lambda_{*} } ( A, B )$ to $\Ext_{ \delta_{*} } ( C, B )$.  We have just proved (1).

Since $K_{*} ( x )$ is an isomorphism, $\image ( \lambda_{*} ) = \image ( \delta_{*} )$.  It is straightforward to show that $K_{j} ( x )$ is an isomorphism from $\ker( \lambda_{j} )$ onto $\ker ( \delta_{j} )$ for $j = 0 ,1$.  Therefore, $[ K_{*} ( x ) ]$ is the isomorphism induced by $K_{*} ( x)$.  This proves (2).

We now prove (3).  Let $z$ be in $\Ext_{ \lambda_{*} } ( C, B )$.  Let
$\e$ be an element of $\Exts ( C, B )$ such that $x_{ C, B } ( \e ) =
z$.  Since $x = \kk ( \alpha )$ for some injective \starhom $\alpha$
from $A$ to $C$, there exist $\alpha \cdot \e$ in $\Exts ( A , B )$
and a homomorphism $( \id_{B} , \eta , \alpha )$ from $\alpha \cdot
\e$ to $\e$.  By \cite{extpurelyinf} Proposition 1.1
\begin{equation*}
x_{ A, B } ( \alpha \cdot \e ) = \kk ( \alpha ) \times x_{ C , B } ( \e ) = x \times z.
\end{equation*}
By the Five Lemma, $( K_{*} ( \id_{B} ) , K_{*} ( \eta ) , K_{*} ( \alpha ) )$ is an isomorphism from $\Ksix ( z )$ onto $\Ksix ( x \times z )$ since $K_{*} ( x )$ and $K_{*} ( \id_{B} )$ are isomorphisms. 

It is clear from the observations made in the previous paragraph and from the definition of $s_{ \delta_{*} }$, $s_{ \lambda_{*} }$, $x \times ( \cdot )$, and $[ K_{*} ( x ) ]$ that the above diagram is commutative.
\end{proof}

R{\o}rdam conjectured that Proposition 3.1 and Theorem 3.2 in
\cite{extpurelyinf} are true for all separable nuclear \cstar-algebras
in \bootstrap. We establish this under the added assumption of stability:

\begin{theorem}\label{T : grphom}
Let $A$ and $B$ be separable nuclear \cstar-algebras in \bootstrap with $B$ stable.  Let $\delta_{*} = ( \delta_{0}, \delta_{1} )$ be an element of $\Hom ( K_{*} ( A ) , K_{  * + 1} ( B ) )$.  
\begin{enumerate}
\item  The map
\begin{equation*}
\ftn{ s_{ \delta_{*}  } = s_{ A, B , \delta_{*} } }{ \Ext_{ \delta_{*} } ( A, B ) } { \Extab ( \ker( \delta_{*} ) , \coker ( \delta_{ * + 1} ) ) }  
\end{equation*}
is a group homomorphism.

\item  If $x$ is in $\Ext ( A, B )$ and if $K _{*} ( x ) = \delta_{*}$, then $s_{ \delta_{*} } ( x ) = \sigma_{ \delta_{*} } ( \Ksix (x ) )$.

\item  If $z$ is in $\Extab ( K_{*} ( A ), K_{*} ( B ) )$, then $s_{ \delta_{*} } ( \epsilon ( z ) ) = \zeta_{ \delta_{*} } ( z )$, where $\epsilon$ is the canonical embedding of $\Extab ( K_{*}( A ) , K_{*} ( B ) )$ into $\Ext ( A, B )$.  
\end{enumerate}
\end{theorem}

\begin{proof}
(2) and (3) are clear from the definition of $s_{ \delta_{*} }$ and $\zeta_{ \delta_{*} }$. 

We now prove (1).  We claim that it is enough to prove (1) for the case that $A$ is a unital separable nuclear purely infinite simple \cstar-algebra in \bootstrap.  Indeed, by the range results in \cite{defKL} and \cite{elliottrordam}, there exists a unital separable nuclear purely infinite simple \cstar-algebra $A_{0}$ in \bootstrap such that $K_{i} ( A )$ is isomorphic to $K_{i} ( A_{0} )$.  Denote this isomorphism by $\lambda_{i}$.  Suppose $A$ is unital.  Then, by Theorem 6.7 in \cite{sepBDF}, there exists an injective $*$-homomorphism $\psi$ from $A$ to $A_{0}$ which induces $\lambda_{*}$.  Suppose $A$ is not unital.  Let $\varepsilon$ be the embedding of $A$ into the unitization of $A$, which we denote by $\widetilde{A}$.  It is easy to find a homomorphism $\widetilde{\lambda}_{i}$ from $K_{i} ( \widetilde{A} )$ to $K_{i} ( A_{0} )$ such that $\widetilde{ \lambda }_{i} \circ K_{i} ( \varepsilon ) = \lambda_{i}$.  Note that $\widetilde{A}$ is a separable unital \cstar-algebra in \bootstrap.  By Theorem 6.7 in \cite{sepBDF}, there exists an injective $*$-homomorphism $\widetilde{ \psi }$ from $\widetilde{A}$ to $A_{0}$ which induces $\widetilde{ \lambda }_{*}$.  Hence, $\psi = \widetilde{ \psi } \circ \varepsilon$ is an injective $*$-homomorphism from $A$ to $A_{0}$ which induces $\lambda_{*}$.  Therefore, in both the unital or the non-unital case, we have an injective $*$-homomorphism $\psi$ which induces an isomorphism from $K_{i} ( A )$ to $K_{i} ( A_{0} )$.  An easy consequence of the Universal Coefficient Theorem \cite{uct} and the Five Lemma shows that $\kk( \psi )$ is a $\kk$-equivalence.  Therefore by Lemma \ref{L : isoextdelta} our claim is true.  

Let $A$ be a unital separable nuclear purely infinite simple \cstar-algebra in \bootstrap.  By the range results of \cite{defKL} and \cite{elliottrordam}, there exist separable nuclear purely infinite simple \cstar-algebras $A_{0}$ and $B_{0}$ in \bootstrap such that $A_{0}$ is unital, $B_{0}$ is stable, and
\begin{align*}
\alpha_{j} : K_{j} ( A_{0} ) \cong \ker ( \delta_{j} : K_{j} ( A ) \to K_{ j+1 } ( B ) ) \\
\beta_{j} : K_{j} ( B_{0} ) \cong \coker ( \delta_{ j+1 } : K_{ j+1 } ( A ) \to K_{j} ( B ) ) 
\end{align*}
for $j = 0 ,1$.  Since $A$ and $A_{0}$ are unital separable nuclear purely infinite simple \cstar-algebras satisfying the Universal Coefficient Theorem, by Theorem 6.7 in \cite{sepBDF} there exists an injective $*$-homomorphism $ \varphi$ from $A_{0}$ to $A$ such that for $j = 0 , 1$ the map $K_{j} ( A_{0} ) \overset{ \alpha_{j} }{ \to } \ker( \delta_{j} ) \hookrightarrow K_{j} ( A )$ is equal to $K_{j} ( \varphi )$.  Choose $b$ in $\kk( B, B_{0} )$ such that for $j = 0, 1$ the map from $K_{j} ( B )$ to $\coker( \delta_{ j + 1 } )$ is equal to $\beta_{j} \circ K_{j}( b )$.  Now, using the same argument as Proposition 3.1 in \cite{extpurelyinf}, we have that the map $s_{ \delta_{*} }$ is a group homomorphism. 
\end{proof}

Replacing Proposition 3.1 in \cite{extpurelyinf} by the above theorem and arguing as in Theorem 3.2 in \cite{extpurelyinf}, we get the following result.

\begin{theorem}\label{T : idealkk}
Let $A$ and $B$ be separable nuclear \cstar-algebras in \bootstrap with $B$ stable.  Suppose $x_{1}$ and $x_{2}$ are elements of $\Ext ( A, B )$.  Then $\Ksix ( x_{1} ) = \Ksix ( x_{2} )$ in $\sixs ( A, B )$ if and only if there exist elements $a$ of $\kk ( A, A )$ and $b$ of $\kk ( B, B )$ with $K_{*} ( a ) = K_{*} ( \id_{A} )$ and $K_{*} ( b ) = K_{*} ( \id_{B} )$ such that $x_{1} \times b = a \times x_{2}$.
\end{theorem}

\section{Classification results}\label{classresult}
We will now use the results of the previous sections to generalize R{\o}rdam's results in \cite{extpurelyinf}.  

Since in the sequel we will be mostly interested in \cstar-algebras that are classified by $(K_{0}(A) , K_{0} ( A )_{+} , K_{1} ( A ) )$, we will not state the Elliott invariant in its full generality.  
 
\begin{definition}\label{D : elliottinv}
For a \cstar-algebra $A$ of real rank zero, the \emph{Elliott invariant} (which we denote by $\Ell ( A )$) consists of the triple 
\begin{equation*}
\Ell ( A ) = ( K_{0} ( A ), K_{0} ( A )_{+}, K_{1} ( A ) ). 
\end{equation*}  
  It is well-known that the canonical embedding of a \cstar-algebra $A$ into its stabilization $A \otimes \K$ induces an isomorphism from $\Ell ( A )$ to $ \Ell ( A \otimes \K )$ (this follows easily from Theorem 6.3.2 and the proof of Proposition 4.3.8 in \cite{bookRordam}).

Suppose $A$ and $B$ are separable nuclear \cstar-algebras in \bootstrap.  Let $x$ be an element of $\kk ( A, B )$.  We say that $x$ induces a homomorphism from $\Ell ( A )$ to $\Ell ( B )$ if $K_{*} ( x )$ is a homomorphism from $\Ell ( A )$ to $\Ell ( B )$.  If, moreover $K_{0} ( x ) ( [ \unit{A} ] ) = [ \unit{B} ]$, then we say $x$ induces a homomorphism from $( \Ell ( A ) , [ \unit{A} ] )$ to $( \Ell ( B ) , [ \unit{B} ] )$.
\end{definition}

\begin{definition}\label{D : classifiable}
We will be interested in classes $\mc{C}$ of separable nuclear unital simple \cstar-algebras in \bootstrap satisfying the following properties:  

\begin{enumerate}[(I)]

\item  Any element of $\mc{C}$ is either purely infinite or stably finite.

\item  $\mc{C}$ is closed under tensoring with $\mbf{M}_{n}$, where $\mbf{M}_{n}$ is the \cstar-algebra of $n$ by $n$ matrices over $\C$.

\item  If $A$ is in $\mc{C}$, then any unital hereditary sub-\cstar-algebra of $A$ is in $\mc{C}$. 

\item  For all $A$ and $B$ in $\mc{C}$ and for all $x$ in $\kk ( A, B )$ which induce an isomorphism from $( \Ell ( A ) , [ \unit{A} ] )$ to $( \Ell ( B ) , [ \unit{B} ] )$, there exists a \stariso $\alpha$ from $A$ to $B$ such that $\kk ( \alpha) = x$.

\end{enumerate}    
\end{definition}

\begin{remark}\label{R : algsinC}

\begin{enumerate}
\item  The class of all unital separable nuclear purely infinite simple \cstar-algebras satisfying the Universal Coefficient Theorem satisfies the properties in Definition \ref{D : classifiable} (see \cite{kirchpureinf} and \cite{phillipspureinf}).
\item  The class of all unital separable nuclear simple
\cstar-algebras satisfying the Universal Coefficient Theorem and with
tracial topological rank zero satisfies
the properties in Definition \ref{D : classifiable} (see Corollary
3.26 in \cite{linnui}). This supersedes earlier work by
Kishimoto-Kumjian (see Corollary 3.13 in \cite{kishkum1}) and by Lin (see
Theorem 1.1
in \cite{morphTAF}). 
\end{enumerate}
\end{remark}


The proof of the following lemma is left to the reader.

\begin{lemma}\label{L : stablehom}
Let $\mc{C}$ be a class of  \cstar-algebras satisfying the properties in Definition \ref{D : classifiable}.  Let $A$ and $B$ be in $\mc{C}$.  Suppose there exists $x$ in $\kk ( A \otimes \K, B \otimes \K )$ such that $x$ induces an isomorphism from $\Ell ( A \otimes \K )$ onto $\Ell ( B \otimes \K )$ and $K_{0} ( x ) (  [ \unit{A} \otimes e_{11} ] ) = [ \unit{B} \otimes e_{11} ]$.  Then there exists a \stariso $\alpha$ from $A \otimes \K$ onto $B \otimes \K$ such that $\kk ( \alpha ) = x$.
\end{lemma}


\begin{lemma}\label{L : reduction}
Let $A_{1} , A_{2} ,B_{1}$, and $B_{2}$ be unital separable nuclear $C\sp*$-algebras and let 
\begin{equation*}
\e : \ 0 \to B_{1} \otimes \K \to E_{1} \to A_{1} \otimes \K \to 0
\end{equation*}
be an essential extension.  Let $\alpha_{*}$ from $\Ell( A_{1} \otimes \K )$ to $\Ell( A_{2} \otimes \K )$ and $\beta_{*}$ from $\Ell( B_{1} \otimes \K )$ to $\Ell( B_{2} \otimes \K )$ be isomorphisms.  Suppose there exist a norm-full projection $p$ in $\mbf{M}_{n} ( A_{1} )$ and a norm-full projection $q$ in $\mbf{M}_{r} ( B_{1} )$ such that $\alpha_{0} ( [ p ] ) = [ \unit{A_{2}} \otimes e_{11} ] $, and $\beta_{0} ( [ q ] ) = [ \unit{B_{2}} \otimes e_{11} ]$.

Then there exist \starisos $\varphi$ from $p \mbf{M}_{n} ( A_{1} ) p \otimes \K $ to $A_{1} \otimes \K$ and $\psi$ from $q \mbf{M}_{r} ( B_{1} ) q \otimes \K $ to $B_{1} \otimes \K$ such that $\varphi \cdot \e$ is isomorphic to $\e$ via the isomorphism $( \id_{ B_{1} \otimes \K }, \id_{E_{1}}, \varphi)$ with $( \alpha_{0} \circ K_{0}( \varphi ) ) ( [ p \otimes e_{11} ] ) = [ \unit{A_{2}} \otimes e_{11} ]$ and $\e$ is isomorphic to $\e \cdot \psi^{-1}$ via the isomorphism $( \psi^{-1} , \id_{ E_{1} } , \id_{A_{1}} )$ with $( \beta_{0} \circ K_{0} ( \psi ) ) ( [ q \otimes e_{11} ] ) = [ \unit{B_{2}} \otimes e_{11} ]$. 

Moreover, $\e$ is isomorphic to $\varphi \cdot \e \cdot \psi^{-1}$ via the isomorphism $( \psi^{-1} , \id_{E_{1}} , \varphi )$.
\end{lemma}

\begin{proof}
By Lemma \ref{L : brown}, there exists a \stariso $\varphi$ from $p ( A_{1} \otimes \K ) p \otimes \K$ to $A_{1} \otimes\K$ such that $[ \varphi ( p \otimes e_{11} ) ] = [ p ]$.  By the definition of $\varphi \cdot \e$, we have that $\varphi \cdot \e$ is isomorphic to $\e$ via the isomorphism $( \id_{ B_{1} \otimes \K }, \id_{E_{1}}, \varphi )$.  Also note that $( \alpha_{0} \circ K_{0}( \varphi ) ) ( [ p \otimes e_{11} ] ) = \alpha_{0} ( [ p ] ) = [ \unit{ A_{2} } \otimes e_{11} ]$.

Using Lemma \ref{L : brown} again, there exists a \stariso $ \psi$ from $q ( B_{1} \otimes \K ) q \otimes \K$ to $B_{1} \otimes \K$ such that $[ \psi ( q \otimes e_{11} ) ] = [ q ]$.  By the definition of $\e \cdot \psi^{-1}$, we have that $\e$ is isomorphic to $\e \cdot \psi^{-1}$ via the isomorphism $(\psi^{-1}, \id_{E_{1}}, \id_{ A_{1} \otimes \K } )$.  Note that $( \beta_{0} \circ K_{0} ( \psi ) ) ( [ q \otimes e_{11} ] ) = \beta_{0} ( [ q ] ) = [ \unit{B_{2}} \otimes e_{11} ]$.
    
Note that the composition of $( \id_{ B_{1} \otimes \K }, \id_{E_{1}}, \varphi )$ with $(\psi^{-1}, \id_{E_{1}}, \id_{ A_{1} \otimes \K } )$ gives an isomorphism $( \psi^{-1} , \id_{E_{1}} , \varphi )$ from $\e$ onto $\varphi \cdot \e \cdot \psi^{-1}$. 
\end{proof}

The next lemma is well-known and we omit the proof.  

\begin{lemma}\label{L : uniteqext}
Let $\e_{1}$ and $\e_{2}$ be in $\Exts ( A , B )$ and let $\tau_{1}$ and $\tau_{2}$ be the Busby invariant of $\e_{1}$ and $\e_{2}$ respectively.  If $\tau_{1}$ is unitarily equivalent to $\tau_{2}$ with implementing unitary coming from the multiplier algebra of $B$, then $\e_{1}$ is isomorphic to $\e_{2}$. 
\end{lemma}

A key component used by R{\o}rdam in \cite{extpurelyinf} was
Kirchberg's absorption theorem.  Elliott and Kucerovsky in
\cite{elliottkuc} give a criterion for when extensions are absorbing and call such extensions \emph{purely large}.  By Kirchberg's theorem, every
essential extension of separable nuclear \cstar-algebras by stable
purely infinite simple \cstar-algebras is purely large. Kucerovsky
and Ng (see \cite{NgCFP} and \cite{KucNgCFPdef}) proved that for 
\cstar-algebras satisfying the corona factorization property, any full
and essential extension is purely large, and hence absorbing.  Properties similar to the corona
factorization property were also studied by Lin \cite{fullext}.

\begin{definition}\label{D : CFP}
Let $B$ be a separable stable \cstar-algebra.  Then $B$ is said to have the \emph{corona factorization property} if every norm-full projection in $\multialg{B}$ is Murray-von Neumann equivalent to $\unit{ \multialg{B} }$.
\end{definition}

The following key results are due to Kucerovsky
and Ng (see \cite{NgCFP} and \cite{KucNgCFPdef}):

\begin{theorem}
Let $A$ be a unital separable simple \cstar-algebra.

\begin{enumerate}

\item If $A$ is exact, $A$ has real rank zero and stable rank one, and $K_{0} ( A )$ is weakly unperforated, then $A \otimes \K$ has the corona factorization property.

\item If $A$ is purely infinite, then $A \otimes \K$ has the corona factorization property.

\end{enumerate}
\end{theorem}

The following theorem is one of two main results in this paper.  Using terminology introduced by Elliott in \cite{elliottclass}, the next result shows that the six term exact sequence together with certain positive cones is a classification functor for certain essential extensions of simple strongly classifiable \cstar-algebras.   

\begin{theorem}\label{T : classification}
Let $\mc{C}_{I}$ and $\mc{C}_{Q}$ be classes of \cstar-algebras satisfying the properties of Definition \ref{D : classifiable}.  Let $A_{1}$ and $A_{2}$ be in $\mc{C}_{Q}$ and let $B_{1}$ and $B_{2}$ be in $\mc{C}_{I}$ with $B_{1} \otimes \K$ and $B_{2} \otimes \K$ satisfying the corona factorization property.  Let 
\begin{align*}
\e_{1} :  \quad \quad & 0 \to B_{1} \otimes \K \to E_{1} \to A_{1} \otimes \K \to 0 \\
\e_{2} :  \quad \quad & 0 \to B_{2} \otimes \K \to E_{2} \to A_{2} \otimes \K \to 0
\end{align*}
be \textbf{essential} and \textbf{full} extensions.  Then the following are
equivalent: 

\begin{enumerate}

\item  $E_{1}$ is isomorphic to $E_{2}$. 

\item   $\e_{1}$ is isomorphic to $\e_{2}$.

\item  There exists an isomorphism $(\beta_{*} , \eta_{*} , \alpha_{*} )$ from $\ksix( \e_{1} )$ to $\ksix ( \e_{2} )$ such that $\beta_{*}$ is an isomorphism from $\Ell ( B_{1} \otimes \K )$ onto $\Ell ( B_{2} \otimes \K )$ and $\alpha_{*}$ is an isomorphism from $\Ell ( A_{1} \otimes \K )$ onto $\Ell ( A_{2} \otimes \K )$.    
\end{enumerate}
\end{theorem}

\begin{proof}
Since $A_{1}$, $A_{2}$, $B_{1}$, and $B_{2}$ are simple
\cstar-algebras, by \cite{extpurelyinf} Proposition 1.2 $E_{1}$ is
isomorphic to $E_{2}$ if and only if $\e_{1}$ is isomorphic to
$\e_{2}$.  It is clear that an isomorphism from $\e_{1}$ onto $\e_{2}$
induces an isomorphism $( \beta_{*} , \eta_{*} , \alpha_{*} )$ from
$\ksix ( \e_{1} )$ onto $\ksix ( \e_{2} )$ such that $\beta_{*}$ is an
isomorphism from $\Ell ( B_{1} \otimes \K )$ onto $\Ell ( B_{2}
\otimes \K )$ and $\alpha_{*}$ is an isomorphism from $\Ell ( A_{1}
\otimes \K )$ onto $\Ell ( A_{2} \otimes \K )$.   

So we only need to prove (3) implies (2).  Using the fact that the canonical embedding of $A_{i}$ into $A_{i} \otimes \K$ induces an isomorphism between $K_{j} ( A_{i} )$ and $K_{j} ( A \otimes \K )$ and since $A_{i}$ is simple,  by Lemma \ref{L : reduction} we may assume $\beta_{0} ( [ \unit{ B_{1} } \otimes e_{11} ] ) = [ \unit{ B_{2} } \otimes e_{11} ]$ and $\alpha_{0} ( [ \unit{ A_{1} } \otimes e_{11} ] ) = [ \unit{ A_{2} } \otimes e_{11} ]$.  Hence, by Lemma \ref{L : stablehom} and the Universal Coefficient Theorem, there exist \starisos $ \beta$ from $B_{1} \otimes \K$ to $B_{2} \otimes \K$ and $\alpha$ from $A_{1} \otimes \K$ to $A_{2} \otimes \K$ such that $K_{*} ( \beta ) = \beta_{*}$ and $K_{*} ( \alpha ) = \alpha_{*}$.   

By \cite{extpurelyinf} Proposition 1.2, $\e_{1}$ is isomorphic to $\e_{1} \cdot \beta$ and $\e_{2}$ is isomorphic to $\alpha \cdot \e_{2}$.  It is straightforward to check that $( K_{*} ( \id_{ B_{2} \otimes \K } ) , \eta_{*} , K_{*} ( \id_{ A_{1} \otimes \K } ) )$ gives a congruence between $\ksix ( \e_{1} \cdot \beta )$ and $\ksix ( \alpha \cdot \e_{2} )$.  Therefore, by Proposition 2.1 in \cite{extpurelyinf} 
\begin{equation*}
\Ksix ( x_{ A_{1} \otimes \K , B_{2} \otimes \K } ( \e_{1} \cdot \beta ) ) = \Ksix ( x_{ A_{1} \otimes \K , B_{2} \otimes \K } ( \alpha \cdot \e_{2} ) ).
\end{equation*}  

Let $x_{j} = x_{ A_{j} \otimes \K , B_{j} \otimes \K } ( \e_{j} )$ for
$j = 1 ,2$.  By \cite{extpurelyinf} Proposition 1.1,
\begin{align*}
\Ksix ( x_{1} \times \kk ( \beta ) ) &= \Ksix (  x_{ A_{1} \otimes \K , B_{2} \otimes \K } ( \e_{1} \cdot \beta ) ) \\
                                    &=  \Ksix ( x_{ A_{1} \otimes \K , B_{2} \otimes \K } ( \alpha \cdot \e_{2} ) ) \\
                                   &= \Ksix ( \kk ( \alpha ) \times x_{2} ). 
\end{align*}      
By Theorem \ref{T : idealkk}, there exist invertible elements $a$ of $\kk ( A_{1} \otimes \K , A_{1} \otimes \K )$ and $b$ of $\kk ( B_{2} \otimes \K , B_{2} \otimes \K )$ such that
\begin{enumerate}
\item $K_{*} ( a ) = K_{*} ( \id_{ A_{1} \otimes \K } )$ and $K_{*} ( b ) = K_{*} ( \id_{ B_{2} \otimes \K } )$ and 
\item $x_{1} \times \kk ( \beta ) \times b = a \times \kk ( \alpha ) \times x_{2}$. 
\end{enumerate} 
Since $A_{1}$ is in $\mc{C}_{Q}$ and $B_{2}$ is in $\mc{C}_{I}$, by
Lemma \ref{L : stablehom} there exist $*$-automorphisms $\rho$ on
$A_{1} \otimes \K$ and $\gamma$ on  $B_{2} \otimes \K$ such that $\kk ( \rho ) = a$ and $\kk ( \gamma )
= b$. 

Using \cite{extpurelyinf} Proposition 1.2 once again, $\e_{1} \cdot
\beta$ is isomorphic to $\e_{1} \cdot \beta \cdot \gamma$ and $\alpha
\cdot \e_{2}$ is isomorphic to $\rho \cdot \alpha \cdot \e_{2}$.  By
\cite{extpurelyinf} Proposition 1.1,
\begin{alignat*}{10}
x_{ A_{1} \otimes \K , B_{2} \otimes \K } ( \e_{1} \cdot \beta \cdot \gamma ) &= x_{1} \times \kk ( \beta ) \times \kk ( \gamma ) & &=x_{1} \times \kk ( \beta ) \times b \\
                                                                              &= a \times \kk ( \alpha ) \times x_{2} &&=\kk ( \rho ) \times \kk ( \alpha ) \times x_{2} \\
                                                                              &= x_{ A_{1} \otimes \K , B_{2} \otimes \K } ( \rho \cdot \alpha \cdot \e_{2} ). & &                                      
\end{alignat*}
  By assumption, $B_2\otimes\mathcal{K}$ satisfies the corona
  factorization property and $\mathfrak{e}_1\cdot\beta\cdot\gamma$ and
  $\rho\cdot\alpha\cdot\mathfrak{e}_2$ are full extensions. 
  The above equation shows that 
  $\mathfrak{e}_1\cdot\beta\cdot\gamma$ and
  $\rho\cdot\alpha\cdot\mathfrak{e}_2$ give the same element of
  $\mathrm{Ext}(A_1\otimes\mathcal{K},B_2\otimes\mathcal{K})$, so by
  Theorem 3.2(2) in \cite{NgCFP} (see also Corollary 1.9 in \cite{KucNgCFPdef}) these
  extensions are unitarily equivalent with the implementing unitary
  coming from the multiplier algebra of $B_2\otimes\mathcal{K}$. 
  So by Lemma \ref{L : uniteqext}, 
  $\mathfrak{e}_1\cdot\beta\cdot\gamma$ is isomorphic to 
  $\rho\cdot\alpha\cdot\mathfrak{e}_2$. Hence 
  $\mathfrak{e}_1$ is isomorphic to $\mathfrak{e}_2$.
\end{proof}

We now extend our results to the case of an ideal which is non-simple
under the added assumption that the ideal is $AF$ but the quotient is
not. We first need:

\begin{lemma}\label{L : afCFP}
Let $A$ be a unital $AF$-algebra.  Then $A \otimes \K$ has the corona factorization property.
\end{lemma}

\begin{proof}
Suppose $p$ is a norm-full projection in $\multialg{ A \otimes \K }$.  Then, by Corollary 3.6 in \cite{fullext}, there exists $z$ in $\multialg{ A \otimes \K }$ such that $z p z^{*} = \unit{ \multialg{ A \otimes \K } }$.  Therefore, $\unit{ \multialg{ A \otimes \K } }$ is Murray-von Neumann equivalent to a sub-projection of $p$.  Since $1_{ \multialg{ A \otimes \K } }$  is a properly infinite projection, $p$ is a properly infinite projection.  By the results of Cuntz in \cite{kthypureinf} and the fact that $K_{0} ( \multialg{ A \otimes \K } ) = 0$, we have that $\unit{ \multialg{ A \otimes \K } }$ is Murray-von Neumann equivalent to $p$. 
\end{proof}

\begin{lemma}\label{L : stableiso}
Let $A$ be a separable stable \cstar-algebra satisfying the corona factorization property.  Let $q$ be a norm-full projection in $\multialg{ A }$.  Then $q A q$ is isomorphic to $A$ and hence $q A q$ is stable.
\end{lemma}

\begin{proof}
Since $q$ is norm-full in $\multialg{ A }$ and since $A$ has the corona factorization property, there exists a partial isometry $v$ in $\multialg{ A }$ such that $v^{*} v = \unit{ \multialg{ A } }$ and $v v^{*} = q$.  Therefore $v$ induces a \stariso from $A$ onto $q A q$.  Since $A$ is stable, $q A q$ is stable.  
\end{proof}

\begin{theorem}\label{T : matsualgs}
Let $\mc{C}$ be a class of 
\cstar-algebras  satisfying the properties of Definition
\ref{D : classifiable} with the further property that it is disjoint
from the class of $AF$ algebras.  Let $A_{1}$ and $A_{2}$ be in $\mc{C}$
and let $B_{1}$ and $B_{2}$ be unital $AF$-algebras.  Suppose 
\begin{align*} 
\e_{1} : \quad \quad 0 \to  B_{1} \otimes \K \overset{ \varphi_{1} }{ \to } E_{1} \overset{ \psi_{1} }{ \to } A_{1} \to 0 \\
\e_{2} : \quad \quad 0 \to  B_{2} \otimes \K \overset{ \varphi_{2} }{\to} E_{2} \overset{ \psi_{2} }{\to}  A_{2} \to 0 	
\end{align*}
are unital essential extensions.  Let $\e_{1}^{s}$ and $\e_{2}^{s}$ be the extensions obtained by tensoring $\e_{1}$ and $\e_{2}$ with the compact operators.  Then the following are equivalent:

\begin{enumerate}
\item  $E_{1} \otimes \K$ is isomorphic to $E_{2} \otimes \K$.

\item  $\e_{1}^{s}$ is isomorphic to $\e_{2}^{s}$.

\item there exists an isomorphism $( \beta_{*} , \eta_{*} ,
\alpha_{*} )$ from $\ksix ( \e_{1} )$ to $\ksix ( \e_{2 } )$ such that
$\beta_{*}$ is an isomorphism from $\Ell ( B_{1})$ to
$\Ell ( B_{2})$ and $\alpha_{*}$ is an isomorphism from
$\Ell ( A_{1} )$ to $\Ell ( A_{2} )$. 
\end{enumerate}
\end{theorem}

\begin{proof}
First we show that (1) implies (2).  Suppose that there exists a
\stariso $\eta$ from $E_{1} \otimes \K$ onto $E_{2} \otimes \K$.  Note
that for $i = 1, 2$, $A_{i} \otimes \K $ is not an $AF$-algebra by
assumption.  Since $[ (\psi_{2} \otimes \id_{
\K }) \circ \eta \circ ( \varphi_{1} \otimes \id_{ \K }  ) ] ( B_{1}
\otimes \K \otimes \K )$ is an ideal of $A_{1} \otimes \K$ and $A_{1}
\otimes \K$ is a simple \cstar-algebra, $[ (\psi_{2} \otimes \id_{ \K
}) \circ \eta \circ ( \varphi_{1} \otimes \id_{ \K }  ) ] ( B_{1}
\otimes \K \otimes \K )$ is either zero or $A_{1} \otimes \K$.  Since
the image of an $AF$-algebra is again an $AF$-algebra, $[ (\psi_{2}
\otimes \id_{ \K }) \circ \eta \circ ( \varphi_{1} \otimes \id_{ \K }
) ] ( B_{1} \otimes \K \otimes \K ) = 0$.  Hence, $\eta$ induces an
isomorphism from $\e_{1}^{s}$ onto $\e_{2}^{s}$. 

Clearly (2) implies both (1) and
\begin{enumerate}
\item[(3')]  there exists an isomorphism $( \beta_{*} , \eta_{*},
\alpha_{*} )$ from $\ksix ( \e_{1}^{s} )$ onto $\ksix ( \e_{2}^{s} )$
such that $\beta_{*}$ is an isomorphism from $\Ell ( B_{1} \otimes \K
\otimes \K )$ onto $\Ell ( B_{2} \otimes \K \otimes \K )$ and
$\alpha_{*}$ is an isomorphism from $\Ell ( A_{1} \otimes \K )$ onto
$\Ell ( A_{2} \otimes \K)$.   
\end{enumerate}
and  as noted in Definition \ref{D : elliottinv}, (3') is equivalent to (3).
We now prove (3') implies (2).   By Lemma \ref{L : reduction}, we may
assume that $\alpha_{0} ( [ \unit{ A_{1} } \otimes e_{11} ] ) = [
\unit{ A_{2} } \otimes e_{11} ]$.  Using strong
classification for $AF$-algebras and for the elements in $\mc{C}$, we get \starisos $\alpha$ from $A_{1} \otimes \K$ to
$A_{2} \otimes \K$ and $\beta$ from $B_{1} \otimes \K \otimes \K$ to
$B_{2} \otimes \K \otimes \K$ such that $K_{*} ( \alpha ) =
\alpha_{*}$ and $K_{*} ( \beta ) = \beta_{*}$.   

By \cite{extpurelyinf} Proposition 1.2, $\e_{1}^{s}$ is isomorphic to $\e_{1}^{s} \cdot \beta$ and $\e_{2}^{s}$ is isomorphic to $\alpha \cdot \e_{2}^{s}$.  It is straightforward to check that $\ksix ( \e_{1}^{s} \cdot \beta )$ is congruent to $\ksix ( \alpha \cdot \e_{2}^{s} )$.  Hence, by Theorem \ref{T : idealkk} there exist invertible elements $a$ of $\kk ( A_{1} \otimes \K, A_{1} \otimes \K )$ and $b$ of $\kk ( B_{2} \otimes \K \otimes \K , B_{2} \otimes \K \otimes \K)$ such that 
\begin{enumerate}[(i)]
\item $K_{*} ( a ) = K_{*} ( \id_{ A_{1} \otimes \K } )$ and $K_{*} ( b ) = K_{*} ( \id_{ B_{2} \otimes \K \otimes \K } )$; and
\item $x_{ A_{1} \otimes \K, B_{2} \otimes \K \otimes \K } ( \e_{1}^{s} \cdot \beta ) \times b = a \times x_{A_{1} \otimes \K, B_{2} \otimes \K \otimes \K}( \alpha \cdot \e_{2}^{s} )$.
\end{enumerate}
By the Universal Coefficient Theorem, $b = \kk ( \id_{ B_{2} \otimes
\K \otimes \K } )$ since $B_{2} \otimes \K \otimes \K$ is an
$AF$-algebra.  By property (3') again combined with  Lemma \ref{L : stablehom} there exists a
$*$-isomorphism $\rho$ from $A_{1} \otimes \K$ to $A_{1} \otimes \K$
such that $\kk ( \rho ) = a$.   

By \cite{extpurelyinf}, Propositions 1.1 and 1.2, $\rho \cdot \alpha \cdot \e_{2}^{s}$ is isomorphic to $\alpha \cdot \e_{2}^{s}$ and 
\begin{align*}
x_{ A_{1} \otimes \K , B_{2} \otimes \K \otimes \K } ( \e_{1}^{s} \cdot \beta ) &= x_{1} \times \kk ( \beta )  \\
                                                                              &= \kk ( \rho ) \times \kk ( \alpha ) \times x_{2} \\
                                                                              &= x_{ A_{1} \otimes \K , B_{2} \otimes \K \otimes \K} ( \rho \cdot \alpha \cdot \e_{2}^{s} ) ,
\end{align*}
where $x_{i} = x_{ A_{i} \otimes \K, B_{i} \otimes \K \otimes \K } ( \e_{i}^{s} )$.

Let $\tau_{1}$ be the Busby invariant of $\e_{1}^{s} \cdot \beta$ and
let $\tau_{2}$ be the Busby invariant of $\rho \cdot \alpha \cdot
\e_{2}^{s}$.  Then, $[ \tau_{1} ] = [ \tau_{2} ] $ in $\Ext ( A_{1}
\otimes \K , B_{2} \otimes \K \otimes \K )$.  Note that $\e_{i}$ is
full by Lemma \ref{littlel}, so by  Proposition \ref{P : stablefull},
so is ${ \e_{i}^{s} }$.  Using this
observation and the fact that $\beta$, $\alpha$, and $\rho$ are
$*$-isomorphisms, it is clear that $\e_{1}^{s} \cdot \beta$ and 
$\rho \cdot \alpha \cdot
\e_{2}^{s}$ are full extensions.

Note that by Lemma \ref{L : afCFP}, $B_{2} \otimes \K \otimes \K$ has
the corona factorization property.  Therefore, by the observations
made in the previous paragraph one can apply Theorem 3.2(2) in
\cite{NgCFP} to get a unitary $u$ in $\multialg{ B_{2} \otimes \K
\otimes \K }$ such that \begin{equation*} \pi ( u ) \tau_{1} (x) \pi (
u )^{*} = \tau_{2} ( x ) \end{equation*} for all $a$ in $A_{1} \otimes
\K$.  Hence, by Lemma \ref{L : uniteqext}, $\e_{1}^{s} \cdot \beta$
and $\rho \circ \alpha \cdot \e_{2}^{s}$ are isomorphic.  Therefore,
$\e_{1}^{s}$ is isomorphic to $\e_{2}^{s}$.  \end{proof}

\begin{remark} Examples of extensions with $AF$ ideals and quotients
which are simple $AD$ algebras of real rank zero are given in
\cite{dadeil} to demonstrate the need of $K$-theory with
coefficients. These examples show that there is no generalization of
the previous theorem to general extensions; one needs to arrange for
fullness for the methods to work. It would be interesting to
investigate if, as suggested by this example, having $K$-theory with
coefficients as part of the invariant could reduce the requirements
on the extension.
\end{remark}

\section{Applications}\label{examples}

Clearly, Theorem \ref{T : classification} applies to essential
extensions of separable nuclear purely infinite simple stable
\cstar-algebras in $\mathcal{N}$ (and gives us the classification
obtained by R\o{}rdam in \cite{extpurelyinf}). We present here two
other examples of classes of special interest, to which our results
apply.

\subsection{Matsumoto algebras}

The results of the previous section apply to a class of \cstar-algebras introduced
in the work by Matsumoto which was investigated in recent work by the
first named author and Carlsen
(\cite{CarlSym},\cite{CarlEilersGraph},\cite{CarlEilersAUG},\cite{CarlEilersKgrpsMat},\cite{CarlEilersOrderedKgrps}).
Indeed, as seen in \cite{CarlSym} we have for each minimal
shift space $\TSS[]$ with a certain technical property $(**)$ introduced in
Definition 3.2 in \cite{CarlEilersKgrpsMat} that the Matsumoto algebra $\mattau[]$ fits in a short
exact sequence of the form
\begin{equation}\label{toke}
\xymatrix{ 0 \ar[r] & \K^n \ar[r] & \mattau[] \ar[r] & C(\TSS) \rtimes_{\sigma} \Z \ar[r] & 0 }
\end{equation}
where $n$ is an integer determined by the structure of the so-called
\emph{special words} of $\TSS$. Clearly the ideal is an $AF$-algebra and by the work of Putnam
\cite{PutnamSTR} the quotient is a unital simple $AT$-algebra with real
rank zero which falls in the class mentioned in Remark \ref{R :
algsinC} (2).  Let us record a couple of consequences of this:

 \begin{corollary}\label{sturmian}
Let $\TSS[\alpha]$ denote the Sturmian shift space associated to the
parameter $\alpha$ in $[0,1] \backslash \Q$ and $\mattau[\alpha]$ the Matsumoto algebra
associated to $\TSS[\alpha]$. If $\alpha$ and $\beta$ are elements of $[0,1] \backslash \Q$, then 
\[
\mattau[\alpha] \otimes \K \cong \mattau[\beta] \otimes\K
\]
if and only if $\Z + \alpha \Z \cong \Z + \beta \Z$ as ordered groups. 
\end{corollary}
\begin{proof}
The extension \eqref{toke} has the six term exact sequence
\[
\xymatrix{ {\Z} \ar[r]^(.35){0} & { \Z + \alpha \Z } \ar@{=}[r] & { \Z + \alpha\Z} \ar[d]\\
{ \Z } \ar@{=}[u] & 0 \ar[l] & {0} \ar[l] }
\]
by Example 5.3 in \cite{CarlEilersOrderedKgrps}. Now apply Theorem \ref{T : matsualgs}.
\end{proof}

The bulk of the work in the papers \cite{CarlEilersGraph}--\cite{CarlEilersOrderedKgrps} is
devoted to the case of shift spaces associated to primitive, aperiodic
substitutions. As a main result, an algorithm is devised to compute
the ordered group $K_0(\mattau[\tau])$ for any such substitution $\tau$,
thus providing new invariants for such dynamical systems up to flow
equivalence (see \cite{ParrySullivan}). The structure result of \cite{CarlSym}
applies in this case as well, and in
fact, as noted in Section 6.4 of \cite{CarlEilersAUG}, the algorithm provides all the
data in the six term exact sequence associated to the extension
\eqref{toke}. This is based on computable objects $\msf{n}_{\tau},\msf{p}_{\tau},\pmb{\msf{A}}_{\tau},\widetilde{ \pmb{ \msf{A} } }_{\tau}$ of which the latter two are square matrices with integer entries. For each such matrix, say $A$ in $\mbf{M}_{n}( \Z )$, we define a group
\[
DG(A)=\lim_{ \to } \left( \xymatrix{ {\Z^n} \ar[r]^-{A} & {\Z^n} \ar[r]^-{A} & \dots } \right)
\]
which, when $A$ has only nonnegative entries, may be considered as an ordered group which will be a dimension group. We get:


\begin{theorem}\label{T : tmcsemain}
 Let $\tau_1$ and $\tau_2$ be basic substitutions, see \cite{CarlEilersKgrpsMat}, over the alphabets $\mfk{a}_1$ and $\mfk{a}_2$, respectively. Then
\[
\mattau[\tau_1] \otimes \K \cong \mattau[\tau_2] \otimes \K
\]
if and only if there exist group isomorphisms $\phi_1,\phi_2,\phi_3$ with $\phi_1$ and $\phi_3$ order isomorphisms, making the diagram
 \begin{equation}\label{orgdiag}
\vcenter{\xymatrix{
{\Z}
\ar[r]^-{\msf{p}_{\tau_1}}
\ar@{=}[d]
&
{\Z^{\msf{n}_{\tau_1}}}
\ar[r]^-{Q_1}
\ar[d]_-{\phi_1}
&
{ DG( \widetilde{ \pmb{ \msf{A} }  }_{\tau_1} ) }
\ar[r]^-{R_1}
\ar[d]_-{\phi_2}
&
{ DG( \pmb{ \msf{ A } } _{\tau_1} ) }
\ar[d]_-{\phi_3}\\
{\Z}
\ar[r]_-{ \msf{p}_{\tau_2} }
&
{ \Z^{ \msf{n}_{\tau_2} } }
\ar[r]_-{Q_2}
&
{ DG( \widetilde{ \pmb{ \msf{ A } }  }_{\tau_2} ) }
\ar[r]_-{R_2}
& { DG( \pmb{ \msf{ A } }_{\tau_2} ) } } }
\end{equation}
commutative.  Here the finite data $\msf{n}_{\tau_i}$ in $\N$,
$\msf{p}_{\tau_i}$ in $\Z^{ \msf{n}_{\tau_i} }$,
$ \pmb{ \msf{ A } }_{\tau_i}$ in $\mbf{M}_{ | \mfk{a}_i |}(\N_0)$, $\widetilde{ \pmb{ \msf{ A } }  }_{\tau_i}$ in
$\mbf{M}_{ |\mfk{a}_i |+\msf{n}_{ \tau_i} } (\Z)$ are as described in \cite{CarlEilersAUG}, the $Q_i$ are defined by the canonical map to the first occurrence of $\Z^{\msf{n}_{\tau_i}}$ in the inductive limit, and $R_i$ are induced by the canonical map from $\Z^{ | \mfk{a}_i | + \msf{n}_{\tau_i} }$ to $\Z^{| \mfk{a}_i |}$. 
\end{theorem}

\begin{proof}
We have already noted above that Theorem \ref{T : matsualgs} applies, proving ``if''. For ``only if'', we use that any \stariso between $\mattau[\tau_1] \otimes \K$ and $\mattau[\tau_2] \otimes \K$ must preserve the ideal in \eqref{toke} and hence induce isomorphisms on the corresponding six term exact sequence which are intertwined by the maps of this sequence as indicated. And since the vectors $\msf{p}_{\tau_i}$ both have all entries positive, the isomorphism $x \mapsto -x$ between $\Z$ and $\Z$ can be ruled out by positivity of $\phi_1$.
\end{proof}

The following reformulation, suggested to us by Takeshi Katsura,
improves the usability of the result:

\begin{corollary}\label{thankstakeshi}
The triple $[K_0(\mattau[\tau]),K_0(\mattau[\tau])_+,\Sigma_\tau]$
with $\Sigma_\tau$ the scale consisting of a multiset in $K_0(\mattau)$ given by
\[
[Q(e_1),\dots,Q(e_{\Ntau})],
\]
($e_i$ the canonical generators of $\Z^{\Ntau}$, and $Q$ the map
defined in Theorem~\ref{orgdiag}),
 is a complete invariant up to stable
isomorphism for the class of Matsumoto algebras $\mattau[\tau]$ associated to  basic substitutions.
\end{corollary}
\begin{proof}
Assume first that a triple $(\phi_1,\phi_2,\phi_3)$ is given as in
\eqref{orgdiag}. By \cite{CarlEilersKgrpsMat}, we may conclude from the fact that $\phi_3$ is an
order isomorphism that the same is true for $\phi_2$. We also note
that $\Ntauo=\Nupsilon$, and that $\phi_1$ must permute the $e_i$ to be
an order isomorphism. Thus, $\phi_2(\Sigma_{\tau_1})=\Sigma_{\tau_2}$.

In the other direction, assume that $\phi_2:K_0(\mattau[\tau_1])\to
K_0(\mattau[{\tau_2}])$ preserves both the positive cone and the
scale. Again, $\Ntauo=\Nupsilon$, and by permuting the generators
according to the identification of $\Sigma_{\tau_1}$ and
$\Sigma_{\tau_2}$
we
get an order isomorphism $\phi_1$ with $Q_2\circ\phi_1=\phi_2\circ Q_1$. Consequently, an isomorphism
$\phi_3$ is induced, and it will be an order isomorphism by
\cite{CarlEilersKgrpsMat}. Finally we see that
\[
\phi_1(\Z\ppvec_{\tau_1})=\phi_1(\ker Q_1)=\ker
Q_2=\Z\ppvec_{\tau_2}
\]
whence we must have $\phi_1(\ppvec_{\tau_1})=\pm \ppvec_{\tau_2}$, and
the negative sign is impossible by positivity.
\end{proof}

Such a classification result puts further emphasis on the
question raised in Section 6.4 in \cite{CarlEilersAUG} of what relation stable
isomorphism of the Matsumoto algebras induces on the shift spaces. We
note here that that relation must be strictly coarser than flow
equivalence:

\begin{example}\label{notflow}
Consider the substitutions
\[
\tau(0) = 10101000 \qquad \tau(1) = 10100
\]
and
\[
\upsilon(0) = 10100100 \qquad \upsilon(1) = 10100
\]
We have that $\mattau[\tau] \otimes \K \cong \mattau[\upsilon] \otimes \K$
although $\TSS[\tau]$ and $\TSS[\upsilon]$ are not flow equivalent.
\end{example}
\begin{proof}
Since both substitutions are chosen to be basic, computations using the algorithm from \cite{CarlEilersGraph} (for
instance using the program \cite{CarlEilersJava}) show that the
invariant reduces to
$[{DG\left(\left[\begin{smallmatrix}5&3\\3&2\end{smallmatrix}\right]\right)},[0]]$
for both substitutions (see Corollary 5.20 in \cite{CarlEilersKgrpsMat}). Hence by Theorem \ref{T : tmcsemain}, the \cstar-algebras $\mattau[\tau]$ and $\mattau[\upsilon]$ are stably isomorphic. However, the \emph{configuration data} (see \cite{CarlEilersGraph}) are different, namely
\[
\xymatrix@R=0.1cm{
{\bullet}\ar@{=}[r]&{\bullet}&&&{\bullet}\ar@{-}[r]&{\bullet}\ar@{-}[l]\\
&&&&{\bullet}\ar@{-}[ur]&{\bullet}\ar@{-}[dl]\\&&&&{\bullet}\ar@{-}[r]&{\bullet,}
}
\]
respectively, and since this is a flow invariant, the shift spaces $\TSS[\tau]$ and $\TSS[\upsilon]$ are not flow equivalent.
\end{proof}

\subsection{Graph algebras}

A completely independent application is presented by the first named
author and Tomforde in a forthcoming paper (\cite{eiltom}) and we sketch
a basic instance of it here. By the work of many hands (see \cite{irbook} and the
references therein) a   \emph{graph $C^*$-algebra} may be associated to any
directed graph (countable, but possibly infinite). When such
$C^*$-algebras are simple, they are always nuclear and in the
bootstrap class $\mathcal N$, and either purely infinite or $AF$.
They are hence, by appealing to either \cite{kirchpureinf} or \cite{af},
classifiable by the Elliott invariant. Our first main result Theorem \ref{T :
classification} applies to prove the following:

\begin{theorem} \label{graph-algebras-graph}
Let $E$ and $E'$ be unital graph algebras with exactly one 
nontrivial ideal $B$ and $B'$, respectively. Then $E\otimes\mathcal{K}\cong 
E'\otimes\mathcal{K}$ if and only if there exists an isomorphism 
$(\beta_*,\eta_*,\alpha_*)$ between the six term exact sequences associated 
with $E$ and $E'$ such that $\alpha_0$ from $K_0(E/B)$ to $K_0(E'/B')$ and $\beta_0$ 
from $K_0(B)$ to $K_0(B')$ are order isomorphisms.
\end{theorem}
\begin{skproof}
Known structure results for graph $C^*$-algebras establish that all of
$B,B',E/B$ and $E'/B'$ are themselves graph $C^*$-algebras, but to
invoke Theorem \ref{T :
classification} we furthermore need to know that $B$ and $B'$ are
stable and of the form $J\otimes \K$ for $J$ a unital graph
algebra. This is a nontrivial result which is established in \cite{eiltom}.

With this  we
can choose as $\mathcal C$ in Theorem \ref{T :
classification} the union of the set of unital, simple, separable, nuclear and
purely infinite algebras with UCT
and the unital simple $AF$-algebras. Then it is easy to check that
properties (I)-(IV) are satisfied, as is the corona factorization property.
\end{skproof}

\begin{example}\label{fourex}
Consider the three graphs
\begin{center}
\includegraphics[width=1.5cm]{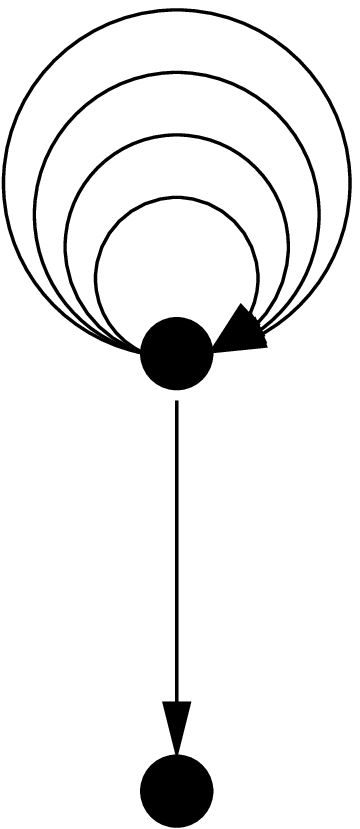}\qquad\includegraphics[width=1.5cm]{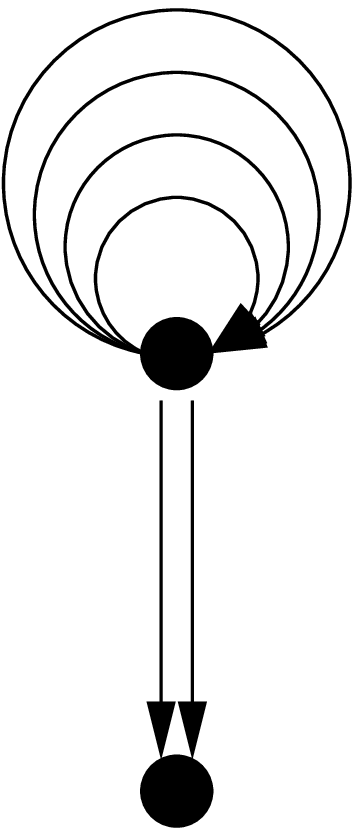}\qquad\includegraphics[width=1.5cm]{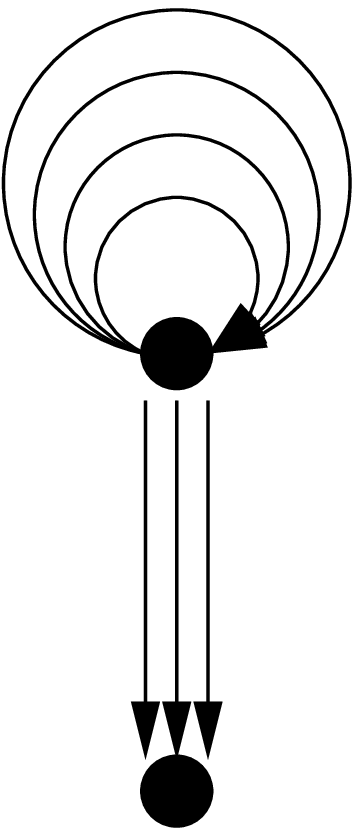}
\end{center}
which all define graph algebras with precisely one ideal, and with
vanishing $K_1$ everywhere. The remaining part of the six term exact
sequence is, up to equivalence in $\operatorname{Ext}(\Z/3,\Z)$,
\[
\xymatrix{{\Z}\ar[r]^-{3}&{\Z}\ar[r]^-{1}&{\Z/3}}\qquad
\xymatrix{{\Z}\ar[r]^-{3}&{\Z}\ar[r]^-{2}&{\Z/3}}\qquad
\xymatrix{{\Z}\ar[r]&{\Z\oplus \Z/3}\ar[r]&{\Z/3}}
\]
respectively. Hence the graph algebras corresponding to the two first
graphs are stably isomorphic to each other, but not to the one
associated to the latter.
\end{example}

This example confirms a special case of a conjecture by Tomforde which we
shall discuss in \cite{eiltom}. Note also that
with substantially more work, Theorem \ref{graph-algebras-graph} above is generalized to
the nonunital case in \cite{eiltom}. An application of
Theorem \ref{T : matsualgs} is also presented there.

\section{Concluding discussion}

\subsection{Extended invariants}

As may be seen by invoking more sophisticated invariants, our assumption of fullness of the extensions considered in
Theorem~\ref{T : classification} is necessary. As is to be expected,
one cannot ignore the order on $K_0$ of the middle algebra in general,
but also invariants such as the tracial simplex or $K$-theory with
coefficients are useful in a non-full context.

It is an interesting and, at present, open question if it is possible
to classify all real rank zero extensions of simple classifiable
$C^*$-algebras by the six-term exact sequence of total $K$-theory.

\subsection{Improved classification results}

It is often of great importance in applications of classification
results to know that the given isomorphism of $K$-groups lifts to a
$*$-isomorphism, or to have direct classification by involving a scale
or a class of the unit in $K$-theory.

The authors in \cite{eilres} and \cite{resrui} have resolved this
question in the case considered by R\o rdam, but the methods used
there do not readily extend to the generality of the present paper.

\subsection{Acknowledgement}

The first author was supported by the EU-Network Quantum Spaces and Noncommutative Geometry (HPRN-CT-2002-00280).

The second author would like to thank the Fields Institute for their hospitality.  Also the second author is grateful for the financial support from the Valdemar Andersen's Travel Scholarship, University of Copenhagen, and the Faroese Research Council.

The second and third author thank Professor George A. Elliott and the participants of the operator algebra seminar at the Fields Institute for many good discussions.

\end{document}